\documentclass[titlepage,twoside,12pt]{article}
\usepackage{amssymb}
\usepackage{amsfonts}
\begin{document}

\pagenumbering{arabic}
\setcounter{page}{1}

\pagenumbering{arabic}

{\LARGE \bf On the Possible Monoid Structures \\ \\
of the Natural Numbers $\mathbb{N}$, \\ \\
or Finding All Associative Binary \\ \\
Operations on $\mathbb{N}$} \\ \\

{\bf Elem\'{e}r E Rosinger \\ Department of Mathematics \\ University of Pretoria \\
Pretoria, 0002 South Africa \\
e-mail : eerosinger@hotmail.com} \\ \\

{\bf Abstract} \\

A certain analysis of all possible {\it associative} binary operations on $\mathbb{N}$ is
presented. This is equivalent with an analysis of all possible {\it monoid} structures on
$\mathbb{N}$. Several results and a conjecture in this regard are given. \\ \\

{\bf 1. Introduction} \\

Let $\mathbb{N} = \{ 1, 2, 3, ~.~.~.~ \}$, and let \\

(1.1)~~~ $ {\cal B} ( \mathbb{N} ) $ \\

denote the set of all {\it binary} operations $f : \mathbb{N} \times \mathbb{N}
\longrightarrow \mathbb{N}$, while by \\

(1.2)~~~ $ {\cal A} ( \mathbb{N} ),~~~ {\cal C} ( \mathbb{N} ) $ \\

we denote the set of all those binary operations $f \in {\cal B} ( \mathbb{N} )$ which are
{\it associative}, respectively, {\it commutative}. \\

By {\it monoid} on $\mathbb{N}$ we mean any structure $( \mathbb{N}, f )$, where $f \in
{\cal A} ( \mathbb{N} )$, thus $f$ is an {\it associative} binary operation on $\mathbb{N}$
which need not necessarily be commutative as well. \\

Obviously $( \mathbb{N}, + )$ and $( \mathbb{N}, . )$ are commutative monoids, the latter also
with neutral element $e = 1$, where "+" and "." denote respectively the usual addition and
multiplication. \\

We define the {\it iterate} of binary operations on $\mathbb{N}$, by the mapping \\

(1.3)~~~ $ {\cal B} ( \mathbb{N} ) \ni f \longmapsto
                  \widetilde f \in {\cal B} ( \mathbb{N} ) $ \\

where \\

(1.4)~~~ $ \widetilde f ( a, b ) ~=~
                       \begin{array}{|l}
                          ~~ a ~~~\mbox{if}~~~ b = 1 \\ \\
                          ~~ f ( \widetilde f ( a, b - 1 ), a ) ~~~\mbox{if}~~~ b \geq 2
                        \end{array} $ \\

It follows that for \\

(1.5)~~~ $ f ( a, b ) = a + b,~~~ a, b \in \mathbb{N} $ \\

we have \\

(1.6)~~~ $ \widetilde f ( a, b ) = a . b,~~~ a, b \in \mathbb{N} $ \\

Further, if now \\

(1.7)~~~ $ f ( a, b ) = a . b,~~~ a, b \in \mathbb{N} $ \\

then \\

(1.8)~~~ $ \widetilde f ( a, b ) = a^b,~~~ a, b \in \mathbb{N} $ \\

With respect to the above iterates $\widetilde f$ of addition and multiplication we can note
the following. \\
First, $\widetilde f$ in (1.6) is both {\it associative} and {\it commutative}. \\
On the other hand, $\widetilde f$ in (1.8) is {\it neither} commutative, {\it nor}
associative. \\
Yet both $\widetilde f$ in (1.6) and (1.8) are {\it distributive} with respect to the
corresponding $f$ in (1.5), respectively (1.7), namely \\

(1.9)~~~ $ \widetilde f ( f ( a, b ), c ) ~=~
         f ( \widetilde f ( a, c ), \widetilde f ( b, c ) ),~~~ a, b, c \in \mathbb{N} $ \\

In Rosinger, the following two {\it uniqueness} properties of the usual {\it addition} $+$ of
natural numbers in $\mathbb{N}$ were proven \\

{\bf Theorem 1.1} \\

Given $f \in {\cal B} ( \mathbb{N} )$. If $f$ is associative and $\widetilde f$ is
commutative, then $f = +$. \\

{\bf Theorem 1.2} \\

Given $f \in {\cal B} ( \mathbb{N} )$. If $f$ is associative and right regular, and
$\widetilde f$ is associative, then $f = +$. \\

Here we used the \\

{\bf Definition 1.1} \\

A binary operation $f \in {\cal B} ( \mathbb{N} )$ is called {\it right regular}, if and
only
if for every $b, c \in \mathbb{N}$, we have \\

(1.10)~~~ $ \left (~~ f ( a, b ) = f ( a, c ),~~~\mbox{for}~~ a \in \mathbb{N} ~~\right )
                ~~~\Longrightarrow~~~ b = c $ \\ \\

{\bf 2. An Open Problem} \\

In view of the role played by associativity in the above Theorems 1.1 and 1.2, one may be
interested in characterizing the set ${\cal A} ( \mathbb{N} )$ of all {\it associative} binary
operations on $\mathbb{N}$. \\
This is equivalent with the characterization of all possible {\it monoid} structures on
$\mathbb{N}$. \\ \\

{\bf 3. ${\cal A} ( \mathbb{N} )$ Is an Infinite Set } \\

Let us start exploring the {\it size} of ${\cal A} ( \mathbb{N} )$. For every {\it bijection}
$~\omega : \mathbb{N} \longrightarrow \mathbb{N}$, we can define the mapping \\

(3.1)~~~ $ {\cal B} ( \mathbb{N} ) \ni f ~\longmapsto~
                             f_\omega \in {\cal B} ( \mathbb{N} )$ \\

by \\

(3.2)~~~ $ f_\omega ( a, b ) ~=~
           \omega ( f ( \omega^{-1} ( a ), \omega^{-1} ( b ) ) ),~~~ a, b \in \mathbb{N} $ \\

It is easy to see that we have \\

{\bf Proposition 3.1} \\

If $~\omega : \mathbb{N} \longrightarrow \mathbb{N}$ is any bijection, then \\

(3.3)~~~ $ f \in {\cal A} ( \mathbb{N} ) ~~\Longrightarrow~~
                                     f_\omega \in {\cal A} ( \mathbb{N} ) $ \\

(3.4)~~~ $ f \in {\cal C} ( \mathbb{N} ) ~~\Longrightarrow~~
                                     f_\omega \in {\cal C} ( \mathbb{N} ) $ \\

{\bf Corollary 3.1} \\

The set ${\cal A} ( \mathbb{N} )$ of associative binary operations on $\mathbb{N}$ is {\it
infinite}. \\

{\bf Proof} \\

The usual addition $+$ and usual multiplication $.$ are obviously in ${\cal A} ( \mathbb{N} )$.
Further, there are infinitely many bijections $~\omega : \mathbb{N} \longrightarrow
\mathbb{N}$. Therefore, there are infinitely many binary operations $(+)_\omega$ and
$(.)_\omega$, each of them associative in view of Proposition 3.1. \\

Indeed, let $\xi, \chi : \mathbb{N} \longrightarrow \mathbb{N}$ be two bijections, and let
us assume that \\

(3.5)~~~ $ (+)_\xi ~=~ (+)_\chi $ \\

Then (3.2) gives $(+)_{(\chi^{-1}\,\circ\,\xi)} = (+)$. But obviously, taking \\

(3.6)~~~ $ \omega^{-1} ~=~ \chi^{-1}\,\circ\,\xi : \mathbb{N} \longrightarrow \mathbb{N} $ \\

we again obtain a bijection, and we have for it \\

(3.7)~~~ $ ( + )_{\omega^{-1}} ~=~ ( + ) $ \\

or in view of (3.2) \\

(3.8)~~~ $ \omega ( a + b ) ~=~ \omega ( a ) + \omega ( b ),~~~ a, b \in \mathbb{N} $ \\

which means that $\omega : \mathbb{N} \longrightarrow \mathbb{N}$ is a {\it homomorphism} of
the {\it monoid} $( \mathbb{N}, + )$. Consequently, if we denote $\omega (1 ) = e \in
\mathbb{N}$, then (3.6) gives \\

$~~~~~~ \begin{array}{l}
          \omega ( 2 ) = \omega ( 1 + 1 ) = \omega ( 1 ) + \omega ( 1 ) = 2 e \\
          \omega ( 3 ) = \omega ( 1 + 2 ) = \omega ( 1 ) + \omega ( 2 ) = e + 2 e = 3 e \\
          \omega ( 4 ) = \omega ( 1 + 3 ) = \omega ( 1 ) + \omega ( 3 ) = e + 3 e = 4 e \\
          .~.~.~.~.~. \\
          \omega ( n ) = n e,~~~ n \in \mathbb{N} \\
          .~.~.~.~.~.
        \end{array} $ \\

However, $\omega$ is a bijection, thus  $\omega$ is {\it surjective}, which means that we must
have $e = 1$, and therefore \\

(3.9)~~~ $ \omega ~=~ id_{\mathbb{N}} $ \\

Now (3.6) implies that $\xi = \chi$. \\

In this way we proved that for any two bijections $\xi, \chi : \mathbb{N} \longrightarrow
\mathbb{N}$, we have \\

(3.10)~~~ $ \xi ~\neq~ \chi ~~~\Longrightarrow~~~ ( + )_\xi ~\neq~ ( + )_\chi $ \\

It follows that $(+)_\omega$, with $\omega : \mathbb{N} \longrightarrow \mathbb{N}$ bijections,
generates infinitely many different elements in ${\cal A} ( \mathbb{N} )$.

\hfill $\Box$ \\

Let us denote by \\

(3.11)~~~ $ bij\, ( \mathbb{N} ) $ \\

the set of all {\it bijections} $~\omega : \mathbb{N} \longrightarrow \mathbb{N}$. Obviously \\

(3.12)~~~ $ ( bij\, ( \mathbb{N} ), \circ ) $ \\

where $\circ$ is the usual composition of functions, is an {\it infinite noncommutative
group}. \\

During the proof of Corollary 3.1, we also proved \\

{\bf Corollary 3.2} \\

The mappings \\

(3.13)~~~ $ \begin{array}{l}
                 bij\, ( \mathbb{N} ) \times {\cal B} ( \mathbb{N} ) \ni ( \omega, f )
                           ~\longmapsto~ f_\omega \in {\cal B} ( \mathbb{N} ) \\ \\

                 bij\, ( \mathbb{N} ) \times {\cal A} ( \mathbb{N} ) \ni ( \omega, f )
                           ~\longmapsto~ f_\omega \in {\cal A} ( \mathbb{N} ) \\ \\
                 bij\, ( \mathbb{N} ) \times {\cal C} ( \mathbb{N} ) \ni ( \omega, f )
                           ~\longmapsto~ f_\omega \in {\cal C} ( \mathbb{N} )
             \end{array} $ \\

are {\it group actions} on ${\cal B} ( \mathbb{N} )$, ${\cal A} ( \mathbb{N} )$ and ${\cal C}
( \mathbb{N} )$, respectively. Furthermore, the mapping \\

(3.14)~~~ $  bij\, ( \mathbb{N} ) \ni \omega ~\longmapsto~
           ( + )_\omega \in {\cal A} ( \mathbb{N} ) ~\bigcap~ {\cal C} ( \mathbb{N} ) $ \\

is {\it injective}. Also, there exist a {\it unique surjective homomorphism} $~\omega :
\mathbb{N} \longrightarrow \mathbb{N}$ of $( \mathbb{N}, + )$, namely \\

(3.15)~~~ $ \omega ~=~ id_{\mathbb{N}} $ \\

which is therefore a {\it bijection} of $\mathbb{N}$. \\ \\

{\bf 4. Monoids on ${\cal B} ( \mathbb{N} )$ } \\

Let $h \in {\cal B} ( \mathbb{N} )$ be any binary operation on $\mathbb{N}$. Then we can
naturally {\it extend} it to a binary operation $h^*$ on ${\cal B} ( \mathbb{N} )$, that is,
to a mapping \\

(4.1)~~~ $ h^* : {\cal B} ( \mathbb{N} ) \times {\cal B} ( \mathbb{N} )
                                   \longrightarrow {\cal B} ( \mathbb{N} ) $ \\

defined by \\

(4.2)~~~ $ h^* ( f, g ) \in {\cal B} ( \mathbb{N} ),~~~ f, g \in {\cal B} ( \mathbb{N} ) $ \\

where \\

(4.3)~~~ $ ( h^* ( f, g ) ) ( a, b ) ~=~ h ( f ( a, b ), g ( a, b ) ),~~~
                                                         a, b \in \mathbb{N} $ \\

Obviously, if we take an {\it associative} binary operation $h \in {\cal A} ( \mathbb{N} )$,
then $( {\cal B} ( \mathbb{N} ), h^* )$ will be a {\it monoid}. And if $h$ is also {\it
commutative}, then so will be the monoid $( {\cal B} ( \mathbb{N} ), h^* )$ . \\

In particular, we can take $h = +$, that is, the usual addition, or $h = .$, which is the
usual multiplication. \\

Also we note that \\

(4.4)~~~ $ h ( a, b ) ~=~ a + b + a . b,~~~ a, b \in \mathbb{N} $ \\

is an {\it associative and commutative} binary operation on $\mathbb{N}$. Therefore, by the
above procedure, it generates a commutative monoid $ ( {\cal B} ( \mathbb{N} ), h^* )$. \\ \\

{\bf 5. An Extended Open Problem} \\

In view of the above, the problem of characterizing the set ${\cal A} ( \mathbb{N} )$ of {\it
associative} binary operations on $ \mathbb{N}$, or equivalently, all the {\it monoid}
structures on $ \mathbb{N}$, can be included in the larger problem of finding all the {\it
monoid} structures on ${\cal B} ( \mathbb{N} )$. \\ \\

{\bf 6. Reformulation of the Open Problem} \\

It is easy to indicate {\it trivial} associative binary operations form ${\cal A}
( \mathbb{N} )$. For instance, for any given fixed $n \in \mathbb{N}$, the {\it constant}
binary operation \\

(6.1)~~~ $ f ( a, b ) ~=~ n,~~~ a, b \in \mathbb{N} $ \\

is obviously associative, as well as commutative. So are the binary operations \\

(6.2)~~~ $ g ( a, b ) ~=~ \min\, \{ a, b \},~~~
                      h ( a, b ) ~=~ \max\, \{ a, b \},~~~ a, b \in \mathbb{N} $ \\

Two other trivial examples of associative binary operations are \\

(6.3)~~~ $ l ( a, b ) ~=~ a,~~~ r ( a, b ) ~=~ b,~~~ a, b \in \mathbb{N} $ \\

both of which, however, fail to be commutative. \\

In view of the above, it is appropriate to restrict the set ${\cal A} ( \mathbb{N} )$ of
associative binary operations to its subset  \\

(6.4)~~~ $ {\cal A}_{mon} ( \mathbb{N} ) $ \\

given by all those associative binary operations $f \in {\cal A} ( \mathbb{N} )$ which are
{\it strictly increasing, separately in each of their two arguments}. \\

Obviously, both the usual addition $+$ and the usual multiplication $.$, as well as the binary
operation in (4.4) which is their sum, belong to ${\cal A}_{mon} ( \mathbb{N} )$. \\

Thus the open problem can be reformulated once more, namely, to characterize all the binary
operations in ${\cal A}_{mon} ( \mathbb{N} )$. \\

It may be useful to consider the following larger class of associative binary operations \\

(6.5)~~~ $ {\cal A}_{gen} ( \mathbb{N} ) $ \\

made up of all those associative binary operations $f \in {\cal A} ( \mathbb{N} )$ which are
{\it genuinely depending on each of their two arguments}, namely, are such that for every $a,
a^{\,\prime}, b, b^{\,\prime} \in \mathbb{N}$, they satisfy the next two conditions \\

(6.6)~~~ $ a \neq a^{\,\prime} ~~\Longrightarrow~~ f ( a, b ) \neq f ( a^{\,\prime}, b ) $ \\

(6.7)~~~ $ b \neq b^{\,\prime} ~~\Longrightarrow~~ f ( a, b ) \neq f ( a, b^{\,\prime} ) $ \\

Obviously \\

(6.8)~~~ ${\cal A}_{mon} ( \mathbb{N} ) ~\subseteq~ {\cal A}_{gen} ( \mathbb{N} )
                                             ~\subseteq~ {\cal A} ( \mathbb{N} ) $ \\

and then the initial open problem can once again be reformulated by looking for a
characterization of any of the above three sets of associative binary operations on
$\mathbb{N}$. \\ \\

{\bf 7. Further Examples of Associative Binary Operations} \\

Suggested by (4.4), let us consider the binary operations of the form \\

(7.1)~~~ $ f ( a, b ) ~=~ \alpha a + \beta b + \gamma a b,~~~ a, b \in \mathbb{N} $ \\

where $\alpha, \beta, \gamma \in \mathbb{N} \cup \{ 0 \}$ are given fixed numbers. Then it is
easy to see that $f$ is associative, if and only if \\

(7.2)~~~ $ \alpha^2 a + \beta c + \alpha \gamma a c =
              \alpha a + \beta^2 c + \beta \gamma a c,~~~ a, b, c \in \mathbb{N} $ \\

which is equivalent with the system \\

(7.3)~~~ $ \alpha^2 = \alpha,~~~ \beta = \beta^2,~~~ \alpha \gamma = \beta \gamma $ \\

that gives the following possible solutions \\

(7.4)~~~ $ \begin{array}{l}
                Case~ 1 :~~~ \alpha = \beta = 1,~~~ \gamma \in \mathbb{N} \cup \{ 0 \} \\ \\
                Case~ 2 :~~~ \alpha = \beta = 0,~~~ \gamma \in \mathbb{N} \cup \{ 0 \} \\ \\
                Case~ 3 :~~~ \alpha = 0,~~~ \beta = 1,~~~ \gamma = 0 \\ \\
                Case~ 4 :~~~ \alpha = 1,~~~ \beta = 0,~~~ \gamma = 0
           \end{array} $ \\

Case 1 contains (4.4), as well as the usual addition, since $\gamma \in \mathbb{N} \cup
\{ 0 \}$ can be arbitrary. Case 2 contains the usual multiplication, due to the same reason.
Case 3 is again trivial since it gives \\

(7.5)~~~ $ f ( a, b ) ~=~ b,~~~ a, b \in \mathbb{N} $ \\

which is associative, but not commutative. The same holds for Case 4, which results in \\

(7.6)~~~ $ f ( a, b ) ~=~ a,~~~ a, b \in \mathbb{N} $ \\

The interesting fact is that we obtain the {\it infinite} family of {\it associative and
commutative} binary operations \\

(7.7)~~~ $ \begin{array}{l}
                   \mathbb{N} \times \mathbb{N} \ni ( a, b) \longmapsto
                           a + b + \gamma a b \in \mathbb{N} \\ \\
                   \mathbb{N} \times \mathbb{N} \ni ( a, b) \longmapsto
                           \gamma a b \in \mathbb{N}
            \end{array} $ \\

with $\gamma \in \mathbb{N} \cup \{ 0 \}$. And clearly, for $\gamma \in \mathbb{N}$, all these
binary operations are in ${\cal A}_{mon} ( \mathbb{N} )$. \\ \\

{\bf 8. A Case of Limitation on Growth} \\

In view of the above, let us check the associativity of binary operations of the form \\

(8.1)~~~ $ f ( a, b ) ~=~ \lambda a^n b^m,~~~ a, b \in \mathbb{N} $ \\

where $\lambda, n, m \in \mathbb{N}$ are given and fixed. It follows easily that such $f$ is
associative, if and only if \\

(8.2)~~~ $ \lambda^n a^{n^2-n} ~=~ \lambda^m c^{m^2-m},~~~ a, c \in \mathbb{N} $ \\

which obviously implies \\

(8.3)~~~ $ n^2 - n ~=~ m^2 - m ~=~ 0 $ \\

thus \\

(8.4)~~~ $ n ~=~ m ~=~ 1,~~~ \lambda \in \mathbb{N} $ \\

and we are back to the second infinite family of binary operation in (7.7). \\

It follows, therefore, that binary operations of type (8.1) must have a growth {\it limited}
to a {\it quadratic} monomial in $a$ and $b$, in order to be associative. \\ \\

{\bf 9. No General Limitation on Growth} \\

The result in section 8 need {\it not} suggest that associative binary operations cannot grow
faster than quadratic monomials in their arguments. Indeed, let us consider any bijection
$~\omega : \mathbb{N} \longrightarrow \mathbb{N}$, such that \\

(9.1)~~~ $ \omega ( 2 n ) ~=~ ( 2 n )^{2 n },~~~ n \in \mathbb{N} $ \\

Then (3.3) implies that $( + )_\omega$ is an associative binary operation on $\mathbb{N}$,
while in view of (3.2), we obtain for $n \in \mathbb{N}$ \\

(9.2)~~~ $ ( + )_\omega ( ( 2 n )^{2 n }, ( 2 n )^{2 n } ) ~=~
       \omega ( 2 \omega^{-1} ( ( 2 n )^{2 n } ) ) ~=~ \omega ( 4 n ) ~=~ ( 4 n )^{4 n } $ \\

In this way \\

(9.3)~~~ $ \limsup_{~a \in \mathbb{N}}~ \frac{(+)_\omega(a,a)}{a^2} ~\geq~
        \limsup_{~n \in \mathbb{N}}~ \frac{( 4 n )^{4 n }}{( 2 n )^{4 n }} ~=~ \infty $ \\

Similarly, $(.)_\omega$ is an associative binary operation on $\mathbb{N}$, and in view of
(3.2), we obtain for $n \in \mathbb{N}$ \\

(9.4)~~~ $ ( . )_\omega ( ( 2 n )^{2 n }, ( 2 n )^{2 n } ) ~=~
       \omega ( ( \omega^{-1} ( ( 2 n )^{2 n } ) )^2 ) ~=~ \omega ( 4 n^2 ) ~=~
                                                             ( 4 n^2 )^{4 n^2 } $ \\

which gives \\

(9.5)~~~ $ \limsup_{~a \in \mathbb{N}}~ \frac{(.)_\omega(a,a)}{a^2} ~\geq~
  \limsup_{~n \in \mathbb{N}}~ \frac{( 4 n^2 )^{4 n^2 }}{( 2 n )^{4 n }} ~=~ \infty $ \\ \\

{\bf 10. A Further Reformulation of the Open Problem} \\

According to Corollary 3.1, the set ${\cal A} ( \mathbb{N} )$ of all the associative binary
operations on $\mathbb{N}$ is infinite, and in view of (3.3) and section 7, it appears to be
rather large. Indeed, there is a large amount of bijections $\omega : \mathbb{N}
\longrightarrow \mathbb{N}$, and for every associative binary operation $f \in {\cal A}
( \mathbb{N} )$ and every such bijection $\omega$, we obtain an associative binary operation
$f_\omega \in {\cal A} ( \mathbb{N} )$. \\

On the other hand, in view of section 9, given an associative binary operation $f \in {\cal A}
( \mathbb{N} )$, we may not necessarily be interested in all the associated binary operations
$f_\omega$, where $\omega$ ranges over all the bijections of $\mathbb{N}$ onto itself. \\

We are thus led to consider the {\it equivalence relation} $\,\approx\,$ on ${\cal B}
( \mathbb{N} )$, defined as follows. Given $f, g \in {\cal B} ( \mathbb{N} )$, then \\

(10.1)~~~ $ f ~\approx~ g ~~~\Longleftrightarrow~~~
            \left(~ \begin{array}{l}
                      \exists~~ \omega \in bij\, ( \mathbb{N} ) : \\ \\
                      ~~~ g ~=~ f_\omega
                     \end{array} ~\right) $ \\

We note that $~\approx~$ is indeed an equivalence relation on ${\cal B} ( \mathbb{N} )$, in
view of Corollary 3.2. Or more directly, (3.2) gives for any $f \in {\cal B} ( \mathbb{N} )$
and for every two bijections $\omega, \chi : \mathbb{N} \longrightarrow \mathbb{N}$, the
relation \\

(10.2)~~~ $ ( f_\omega )_\chi ~=~ f_{( \chi\, \circ\, \omega )} $ \\

and obviously $\chi\, \circ\, \omega : \mathbb{N} \longrightarrow \mathbb{N}$ is again a
bijection. \\

In this way, we can reformulate the open problem by asking to characterize the quotient set \\

(10.3)~~~ $ {\cal A}_{mon} ( \mathbb{N} ) / \approx $ \\

We further note the obvious property of bijections of $\mathbb{N}$ onto itself, namely \\

(10.4)~~~ $ \begin{array}{l}
                 \forall~~ \omega : \mathbb{N} \longrightarrow \mathbb{N} : \\ \\
                 ~~~ \omega ~~\mbox{increasing bijection}
                                  ~~~\Longrightarrow~~~ \omega ~=~ id_{\mathbb{N}}
             \end{array} $ \\

Thus we can ask whether the following property may hold \\

(10.5)~~~ $ \begin{array}{l}
               \forall~~ f, g \in {\cal A}_{mon} ( \mathbb{N} ) : \\ \\
               ~~~ f ~\approx~ g ~~~\Longrightarrow~~~ f ~=~ g
            \end{array} $ \\

or perhaps, the property \\

(10.6)~~~ $ \begin{array}{l}
               \forall~~ f, g \in {\cal A}_{mon} ( \mathbb{N} ),~
                  \omega : \mathbb{N} \longrightarrow \mathbb{N} ~~\mbox{bijection} : \\ \\
               ~~~ g ~=~ f_\omega ~~~\Longrightarrow~~~ \omega ~~\mbox{increasing}
            \end{array} $ \\

We note that in view of (10.4), we have the implication (10.6) $\Longrightarrow$ (10.5). \\

In case (10.5) or (10.6) holds, it follows that each equivalence class in the quotient in
(10.3) contains only one single element. \\ \\

{\bf 11. The Conjecture} \\

Related to the various formulations of the above open problem, we make the conjecture \\

(11.1)~~~ $ {\cal A}_{mon} ( \mathbb{N} ) ~=~
            \{~ f ~~\mbox{in}~ (7.7), ~~\mbox{for}~ \gamma \in \mathbb{N} ~\}
                                     ~\cup~ \{~ + ~\} $ \\

in other words, ${\cal A}_{mon} ( \mathbb{N} )$ consists of the usual addition $+$, usual
multiplication $.$, and the binary operations \\

(11.2)~~~ $ \mathbb{N} \times \mathbb{N} \ni ( a, b) \longmapsto
                           a + b + \gamma a b \in \mathbb{N} $ \\

with $\gamma \in \mathbb{N}$, as well as \\

(11.3)~~~  $ \mathbb{N} \times \mathbb{N} \ni ( a, b) \longmapsto
                                         \gamma a b \in \mathbb{N} $ \\

where $\gamma \in \mathbb{N},~ \gamma \geq 2$. \\ \\

{\bf 12. Associativity and Submonoids} \\

Given $f \in {\cal B} ( \mathbb{N} )$ and $a \in \mathbb{N}$, we define $f_{a\,\bullet},
f_{\bullet\,a} : \mathbb{N} \longrightarrow \mathbb{N}$, by \\

(12.1)~~~ $ f_{a\,\bullet} ( b ) ~=~ f ( a, b),~~~
                       f_{\bullet\,a} ( b ) ~=~ f ( b, a ),~~~ b \in \mathbb{N} $ \\

The next {\it characterization} of associative binary operations follows immediately \\

{\bf Lemma 12.1} \\

For $f \in {\cal B} ( \mathbb{N} )$ we have \\

(12.2)~~~ $ \begin{array}{l}
                f \in {\cal A} ( \mathbb{N} ) ~~~\Longleftrightarrow~~~
                \left(~~ f_{a\,\bullet} \circ f_{b\,\bullet} ~=~ f_{f(a,b)\,\bullet},~~~
                a, b \in \mathbb{N} ~~\right) ~~~\Longleftrightarrow~~~ \\ \\
                ~~~~~~~~~\Longleftrightarrow~~~
                \left(~~ f_{\bullet\,a} \circ f_{\bullet\,b} ~=~ f_{\bullet\,f(b,a)},~~~
                a, b \in \mathbb{N} ~~\right)
             \end{array} $ \\

where $\circ$ is the usual composition of functions.

\hfill $\Box$ \\

Consequently, for $f \in {\cal B} ( \mathbb{N} )$, let us denote \\

(12.3)~~~ $ \begin{array}{l}
                 {\cal R}_f ~=~ \{~ f_{a\,\bullet} ~~|~~ a \in \mathbb{N} ~\} \\ \\
                 {\cal L}_f ~=~  \{~ f_{\bullet\,a} ~~|~~ a \in \mathbb{N} ~\}
             \end{array} $ \\

Therefore we obtain \\

{\bf Corollary 12.1} \\

If $f \in {\cal A} ( \mathbb{N} )$, then $( {\cal R}_f, \circ )$ and $( {\cal L}_f, \circ )$
are {\it submonoids~} in $( \mathbb{N}^{\,\mathbb{N}}, \circ )$.

\hfill $\Box$ \\

In view of the conjecture in section 11, let us denote by \\

(12.4)~~~ $ mon\, ( \mathbb{N} ) $ \\

the set of all {\it strictly increasing} mappings $f : \mathbb{N} \longrightarrow \mathbb{N}$.
Obviously $( mon\, ( \mathbb{N} ), \circ )$ is a {\it semigropup}, that is, a monoid with
neutral element. \\

It follows that we have \\

{\bf Corollary 12.2} \\

Given $f \in {\cal A}_{mon} ( \mathbb{N} )$, then \\

(12.5)~~~ $ ( {\cal R}_f, \circ ) ~\mbox{and}~ ( {\cal L}_f, \circ )
                ~\mbox{are submonoids in}~ ( mon\, ( \mathbb{N} ), \circ ) $

\hfill $\Box$ \\

The result in Corollary 12.1 leads also to \\

{\bf Corollary 12.3} \\

We have the {\it injective} mappings \\

(12.6)~~~ $ \begin{array}{l}
               {\cal A} ( \mathbb{N} ) \ni f \longmapsto \lambda_f : \mathbb{N}
                           \longrightarrow \mathbb{N}^{\mathbb{N}} \\ \\
               {\cal A} ( \mathbb{N} ) \ni f \longmapsto \rho_f : \mathbb{N}
                           \longrightarrow \mathbb{N}^{\mathbb{N}}
             \end{array} $ \\

where $\lambda_f, \rho_f$ are {\it monoid homomorphism} between $( \mathbb{N}, f )$ and
$( \mathbb{N}^{\mathbb{N}}, \circ )$, defined by \\

(12.7)~~~ $ \lambda_f ( a ) ~=~ f_{a\,\bullet},~~~ \rho_f ( a ) ~=~ f_{\bullet\,a},~~~
                                                                 a \in \mathbb{N} $ \\

{\bf Remark 12.1} \\

A likely consequence of (12.5) appears to be the following {\it dichotomy}. For any given
associative binary operation $f \in {\cal A}_{mon} ( \mathbb{N} )$

\begin{itemize}

\item either the growth of $f_{a\,\bullet}$ and $f_{\bullet\,a}$, with $a \in \mathbb{N}$, is
{\it linear}, as conjectured in section 11,

\item or on the contrary, it is {\it much more fast}.

\end{itemize}

Indeed, let us take for instance $g \in mon\, ( \mathbb{N} )$, given by \\

(12.8)~~~ $ g ( n ) ~=~ \alpha\, n^{\,\beta},~~~ n \in \mathbb{N} $ \\

where $\alpha, \beta \in \mathbb{N}$ are arbitrary but fixed. Then clearly \\

(12.9)~~~ $ ( g \circ g ) ( n ) ~=~ \alpha^{\,\beta + 1}\, n^{\,\beta^{\,2}},~~~
                                                              n \in \mathbb{N} $ \\

thus we have \\

(12.10)~~~ $ g \circ g ~\leq~ C\, g $ \\

for some constant $0 < C < \infty$, if and only if \\

(12.11)~~~ $ \beta ~=~ 1 $ \\

It follows that one possible way to prove the conjecture in section 11 is by showing that the
second alternative in the above dichotomy is not possible. \\

A first question in this regard is the following : is it true that \\

(12.12)~~~ $ \begin{array}{l}
                 \forall~~ g \in mon\, ( \mathbb{N} ) : \\ \\
                 \exists~~ f \in {\cal A}_{mon} ( \mathbb{N} ) : \\ \\
                 ~~~ g \in {\cal R}_f ~\cup~ {\cal L}_f
             \end{array} $ \\

In case the answer is "yes", then the above second alternative is possible, thus the
conjecture in section 11 does {\it not} hold. \\

{\bf Remark 12.2} \\

In view of the possible relevance of arguments based on growth, it may be useful to consider
{\it continuous}, and in fact, {\it differentiable} versions of binary, and in particular,
associative binary operations, like for instance \\

(12.13)~~~ $ f : ( 0, \infty ) \times ( 0, \infty ) ~\longrightarrow~ ( 0, \infty ) $ \\

Let us therefore denote by \\

(12.14)~~~ $ {\cal A} ( 0, \infty ) $ \\

the set of all {\it associative} binary operations $f$ in (12.13). \\

A difficulty arises here from the fact that there is {\it no} simple natural way to relate the
two sets ${\cal A} ( \mathbb{N} )$ and ${\cal A} ( 0, \infty )$, although the domains of their
respective elements are in the obvious relation \\

$~~~~~~ \mathbb{N} \times \mathbb{N} ~\subset~ ( 0, \infty ) \times ( 0, \infty ) $ \\

Indeed, given $f \in {\cal A} ( 0, \infty )$, if we consider its restriction to $\mathbb{N}
\times \mathbb{N}$, then we need not obtain an integer valued function, let alone one in
${\cal A} ( \mathbb{N} )$. Conversely, if we take $f \in {\cal A} ( \mathbb{N} )$, such a
function may not always be extendable to one defined on the whole of $( 0, \infty ) \times
( 0, \infty )$, and which at the same time is associative. \\

Nevertheless, we shall briefly consider the continuous case as well. In this regard, we adapt
accordingly and in an obvious manner the notations in (1.1), (1.2), (6.4), etc. \\

Let us introduce the following classes of binary operations on $( 0, \infty)$, according to
their respective {\it growth} conditions. For $\alpha, \beta \in \mathbb{R}$, we denote by \\

(12.15)~~~ $ {\cal B}^{\,\alpha,\,\beta}_{pol} ( 0, \infty ) $ \\

the set of all binary operations $f \in {\cal B} ( 0, \infty )$ which have the {\it polynomial
growth} property \\

(12.16)~~~ $ \begin{array}{l}
                \forall~~ a, b \in ( 0, \infty ) : \\ \\
                \exists~~ K(a,b), L(a,b), H(a,b) \in ( 0, \infty ) : \\ \\
                \forall~~ c \in ( 0, \infty ),~ c \geq H(a,b) : \\ \\
                ~~~ K(a,b)\, c^\alpha f ( a, b ) ~\leq~ f ( a, c b ) ~\leq~
                                 L(a,b)\, c^\beta f ( a, b)
             \end{array} $ \\

Let us now take any \\

(12.17)~~~ $ f \in {\cal B}^{\,\alpha,\,\beta}_{pol} ( 0, \infty )
                                       \cap {\cal A}_{mon} ( 0, \infty ) $ \\

and $a, b, c, d \in ( 0, \infty )$. \\

In the sequel, we consider that $d$ is sufficiently large. \\

Then we have in view of (12.16) applied to $f(f(a,b),dc)$, the inequalities \\

(12.18)~~~ $ \begin{array}{l}
                     K(f(a,b),c)\,d^\alpha f(f(a,b),c) \leq f(f(a,b),dc) \leq \\ \\
                       ~~~~~~~~~~~~~~~~~~~~~~~~~~~~~ \leq L(f(a,b),c)\,d^\beta f(f(a,b),c)
              \end{array} $ \\

But the associativity of $f$ gives \\

(12.19)~~~ $ f(f(a,b),dc) = f(a,f(b,dc))$ \\

while (12.16) applied to $f(b,dc)$ results in \\

(12.20)~~~ $ K(b,c)\, d^\alpha f(b,c) \leq f(b,dc) \leq L(b,c)\,d^\beta f(b,c) $ \\

Therefore (12.20) and the monotonicity of $f$ yield \\

(12.21)~~~ $ f(a,K(b,c)\, d^\alpha f(b,c)) \leq f(a,f(b,dc)) \leq
                                                 f(a,L(b,c)\,d^\beta f(b,c)) $ \\

And now (12.16) applied to both $f(a,K(b,c)\, d^\alpha f(b,c))$ and \\ $f(a,L(b,c)\,d^\beta
f(b,c))$, together with (12.21), will imply \\

(12.22)~~ $ \begin{array}{l}
              K(a,f(b,c)) ( K(b,c)\,d^\alpha )^\alpha f(a,f(b,c)) \leq f(a,f(b,dc)) \leq \\ \\
                   ~~~~~~~~~~~~~ \leq L(a,f(b,c)) ( L (b,c)\,d^\beta )^\beta f(a,f(b,c))
            \end{array} $ \\

In view of (12.19), we can compare (12.18) and (12.22), and thus obtain \\

(12.23)~~~ $ \begin{array}{l}
                K(a,f(b,c)) ( K(b,c)\,d^\alpha )^\alpha \leq L(f(a,b),c)\,d^\beta \\ \\
                K(f(a,b),c)\,d^\alpha \leq L(a,f(b,c)) ( L (b,c)\,d^\beta )^\beta
             \end{array} $ \\

hence since $d \in ( 0, \infty )$ can be arbitrarily large, it follows that \\

(12.24)~~~ $ \alpha^2 \leq \beta,~~~ \alpha \leq \beta^2 $ \\

Of course, in view of (12.16), we have \\

(12.25)~~~ $ \alpha \leq \beta $ \\


\begin{thebibliography}{99}

\bibitem{} Rosinger, E E : The algebraic uniqueness of the addition of natural numbers.
Aequationes Mathematicae, Vol. 25, 1982, 269-273

\end{thebibliography}
\end{document}